\documentclass[numbook,envcountreset,final]{svjour3}

\usepackage[english]{babel}
\usepackage{graphicx,epsf,amssymb,amsmath,rotating}
\usepackage{subeqnarray,mathabx,upgreek,textgreek}
\usepackage{showlabels}

\spnewtheorem{algorithm}[theorem]{Algorithm}{\bf}{\rm}
\spnewtheorem{procedure}[theorem]{Procedure}{\bf}{\rm}

\smartqed

\def\no{\noindent}
\def\rit{\mathbb{R}}

\def\endprf{\qed}
\def\prf{\no{\it Proof.\enspace}}
\def\personalitem{{\boldmath$\cdot$\;}}
\def\fmodel{\widecheck{f}}
\def\bundleset{\mathcal B}
\def\X{\mathrm X}

\begin{document}

\author{Adam Ouorou}
\title{Proximal bundle algorithms for nonsmooth convex optimization
 via fast gradient smooth methods}
\institute{
Orange Labs Research,
44 avenue de la R\'epublique, 
92300 Chatillon, France. \\
\email{adam.ouorou@orange.com}
}


\maketitle

\begin{abstract}
We propose new  proximal bundle algorithms for minimizing a
nonsmooth convex function. 
These algorithms are derived
from the application of Nesterov fast gradient methods for smooth
convex minimization to the so-called Moreau-Yosida regularization
$F_\mu$ of $f$ w.r.t. some $\mu>0$. 
Since the exact values and
gradients of $F_\mu$ are difficult to evaluate, we use approximate proximal
points thanks to a bundle strategy to get implementable
algorithms. One of these algorithms appears as an implementable version
of a special case of inertial proximal algorithm. We
give their complexity estimates in terms of the original function
values, and report some preliminary numerical results.
\end{abstract}


\section{Introduction}
We consider the problem
\begin{equation}\label{pb}
\min\limits_{x\in\rit^n} f(x),
\end{equation}
where $f$ is a convex (non necessarily differentiable) function. We
assume that the set  $\X^*$ of minimizers of $f$ is nonempty.  It is
well known that this problem can be transformed into a differentiable
convex minimization problem
\begin{equation}\label{MYpb}
\min\limits_{x\in\rit^n} F_\mu(x),
\end{equation}
where $\mu >0$, $\|.\|$ is the usual Euclidean norm, and $F_\mu$ is the 
Moreau-Yosida regularization of $f$ defined by 
\begin{equation}\label{MY}
F_{\mu}(x)=\min\limits_{z\in\rit^n}\left\{f(z)+\dfrac{\mu}{2}\|z-x\|^2\right\}.
\end{equation}
The parameter $\mu$ is usually termed as the
\emph{proximity parameter}. 
The function $F_{\mu}$ is a differentiable convex function defined on the whole space $\rit^n$ and has 
$\mu$-Lipschitzian gradient without any further assumption
\cite{HUL93}, i.e.,
$\|\nabla F_\mu(x)-\nabla F_\mu(y)\|\le\mu\|x-y\|,\;\;x,y\in\rit^n$. 
The unique minimizer in 
\eqref{MY} is called the \emph{proximal point} of $x$ and we denote it
by $p_{\mu}(x)$, i.e. 
\begin{equation}\label{prox-pt}
p_{\mu}(x)=\arg\min\limits_{z\in\rit^n}\left\{f(z)+\dfrac{\mu}{2}\|z-x\|^2\right\}.
\end{equation}
The derivative of $F_{\mu}$ is given by 
\begin{equation}\label{Fgrad}
\nabla F_{\mu}(x)=\mu(x-p_{\mu}(x)),
\end{equation}
and $\nabla F_{\mu}(x)\in \partial f(p_{\mu}(x))$ where $\partial f$
is the subdifferential of $f$.  
Minimizing $f$ and $F_{\mu}$ are equivalent problems, 
in the sense that the minima of the two functions coincide, see [\cite{HUL93}, Theorem XV.4.1.7]. 
Assuming a fast
computation of its gradient (in fact $p_\mu(x)$), an efficient smooth minimization algorithm applied to
$F_{\mu}(x)$ is appealing (attractive). This explains the motivation
of developping quasi-Newton type
algorithms for the minimzation of $F_\mu$, see for
instance \cite{BGLS95,LeSa97,Mif96,ChFu99}.  The proximal
point algorithm \cite{Roc76} for solving \eqref{pb} is as follows.
\begin{algorithm}\label{ppa}
\begin{enumerate}
\item[] 
\hspace*{-0.5cm}{\tt Proximal Point Algorithm (PPA)}
\item[0.]\label{ppa0} Choose $x^0\in\rit^n$ and set $k=0$.  
\item\label{ppa1}Compute $p_\mu(x^k)$.
\item\label{ppa2} If $p_{\mu}(x^k)=x^k$ stop: $x^k$ solves \eqref{pb}.
\item\label{ppa3} $x^{k+1}=p_\mu(x^k)$. Increase $k$ by 1 and loop to Step \ref{ppa1}.
\end{enumerate}
\end{algorithm}
As already observed in the literature, the proximal point algorithm
can be 
regarded as a standard gradient algorithm applied to the
minimization of $F_{\mu}$. The classical gradient descent (CGDA) is one of the simplest
method for smooth convex minimization. It writes $x^{k+1}=x^k-\alpha_k\nabla F_{\mu}(x^k)$ 
for \eqref{MYpb} where $\alpha_k$ is a stepsize, and stops when $\nabla F_{\mu}(x^k)=0$.  
\begin{algorithm}\label{gm}
\begin{enumerate}
\item[] 
\hspace*{-0.5cm}{\tt Classical Gradient Descent Algorithm (CGDA)}
\item[0.]\label{gm0} Choose $x^0\in\rit^n$ and set $k=0$.  
\item\label{gm1}Compute $\nabla F_\mu(x^k)$.
\item\label{gm2} If $\nabla F_{\mu}(x^k)=0$ stop: $x^k$ solves \eqref{pb}.
\item\label{gm3} $x^{k+1}=x^k-\alpha_k\nabla F_{\mu}(x^k)$. Increase $k$ by 1 and loop to Step \ref{gm1}.
\end{enumerate}
\end{algorithm}
There are different strategies of choosing the stepsize 
$\alpha_k$, leading to various versions of CGDA.  Since $\nabla F_{\mu}(x^k)=\mu(x^k-p_{\mu}(x^k))$, the stopping
criterion in this algorithm is exactly the same as in PPA. By setting 
$\alpha_k=\mu^{-1},\;k\ge 0$, CGDA reduces to PPA. 
The differentiability of $F_\mu$ motivates us to investigate alternatives
to classical gradient methods which are simple but not optimal \cite{Nes04}. 
In this paper, we consider fast gradient methods initiated by Nesterov
in \cite{Nes83,Nes04}, which attain 
the optimal oracle complexity for smooth convex
optimization. 
Their remarkable feature is that, as in a classical gradient method, they do not need more than one gradient evaluation at
each iteration. The development of
fast first-order methods for smooth problems is an active area of research
\cite{BeTe09,GoKa13,DGN14,KiFe16,Rud17}, motivated by the need to solve
large scale problems unsuited to second-order methods (so is Problem \eqref{MYpb} as $F_\mu$ is not twice differentiable in general \cite{LeSa97b}). The idea of exploiting these
fast methods for the optimization of nonsmooth convex functions is not
new. There is an increasing interest in the context of computing the
zeros of the sum of a maximally monotone operator, resulting in the
class of of so-called \emph{inertial proximal algorithms}, see for
instance \cite{AlAt01,ACR18,ACR19}
and references therein. In \cite{Gul92}, G\"uler extended the concept of estimate sequences (see
[\cite{Nes04}, Definition 2.2.1]) to the nonsmooth function $f$ from
which a main algorithm is established with a convergence rate 
estimate $O(1/k^2)$. 
This algorithm is conceptual in the sense that it makes use
of the exact solutions of the same type of problems as \eqref{prox-pt} 
for some $x=x^k$. We already pointed out the difficulty to solve these
problems in practice. 
A variant in which approximate proximal points can be used 
has been proposed by G\"uler, according to the following criterion
proposed by Rockafellar in \cite{Roc76} to compute an
approximate proximal point $z^{k+1}$ for a given point $x^k$,
\begin{equation}\label{rock-cond}
\min\left\{\|g\|:\;\;g\in\partial\phi_k(z^{k+1})\right\}\le\frac{\varepsilon_k}{\lambda_k},
\end{equation}
where $\varepsilon_k=O(k^{-\sigma})$ for some $\sigma>0$ and
$\phi_k(z)=f(z)+\frac{1}{2\lambda_k}\|z-x^{k}\|^2,\;\lambda_k>0$. 
The l.h.s in \eqref{rock-cond} expresses the distance of $0$ to the
set $\partial\phi_k(z^{k+1})$. The
criterion \eqref{rock-cond} is not always easy to check in
practice since the l.h.s problem is not tractable most of the time. No numerical experiment has been conducted in \cite{Gul92}
to have an idea on the practical efficiency of the proposed
algorithms. 
Our approach is closed in spirit to that of G\"uler in
\cite{Gul92} 
but our purpose is to propose implementable algorithms for a wide
range of problems of type \eqref{pb}, that make use
of approximate function and gradient values of $F_\mu$. 
To this aim, we simply use a bundle strategy to perform these approximate
computations. 

We use some standard notations throughout the paper. The symbol $\langle .,. \rangle$ denotes the usual 
scalar product while the Euclidean norm is denoted by $\|.\|$. For any
$\epsilon\ge 0$, the $\epsilon$-subdifferential of $f$ at $z$ is 
$\partial_{\epsilon} f(z)=\{g:\;f(x)\ge f(z)+\langle g,\;x-z\rangle -\epsilon\;\forall x\in\rit^n$. 
This set is identical with 
the subdifferential of $f$ at $z$ when $\epsilon=0$.

The paper is organized as follows. In the next section, we present
formally a new class of proximal bundle algorithms. In
Section~\ref{sec-conv}, we give their complexity estimates and analyze
in Section~\ref{sec-err}
the accumulation of errors due to the inexact computation of gradients
of $F_\mu$, and propose practical tolerances for these computations. In Section~\ref{sec-num}, we report some preliminary
computational results obtained with the proposed algorithms and
conclude in Section~\ref{sec-con}. 

\section{New algorithms}
We recall that $F_\mu$ is $\mu$-smooth i.e. $\nabla F_\mu$ is
Lipschitz continuous with constant $\mu$. We also have by definition of $F_\mu$, 
$F_\mu(x)\le f(x)\;\text{for any}\;x\in\rit^n$. 
We will make use of the
following properties relating the two problems \eqref{pb} and
\eqref{MYpb}, see for instance [\cite{HUL93}, Theorem XV.4.1.7].
\begin{proposition}\label{MYprop}
The following statements are equivalent: 
\text{\em a) }$x$ minimizes $f$, \text{\em b) }$x=p_\mu(x)$, \text{\em c) }$\nabla F_\mu(x)=0$, \text{\em d) }$x$ minimizes $F_\mu$, \text{\em e) }$f(x)=f(p_\mu(x))$,
\text{\em f) }$f(x)=F_\mu(x)$.
\end{proposition}
The Moreau-Yosida regularization $F_\mu$ provides a smooth lower approximation
of $f$ which coincides with $f$ at optimality. One can then apply a
fast gradient method to $F_\mu$ in order to get a minimizer of $f$.  Based on the above proposition, a \emph{gap function} may be defined as
$\delta(x)=f(x)-F_\mu(x)\;[\ge 0]$, which gives $\delta(x)=0$ iff $x$ is optimal for \eqref{pb}. 
\subsection{Derivation}
The fast gradient
method developped in \cite{Nes83} for smooth convex  functions, uses the 
sequence of reals
$$
\lambda_0=1,\;\;\lambda_{k+1}=\dfrac{1+\sqrt{1+4\lambda_{k}^2}}{2},\;\;k\ge
0,
$$
which satisfies the following useful relations
\begin{equation}\label{lrel}
\lambda_{k-1}^2=\lambda_{k}(\lambda_{k}-1),\;k\ge 1,
\end{equation}
and 
\begin{equation}\label{lrel2}
\lambda_k^2=\sum\limits_{i=0}^k\lambda_i,\;\;\lambda_k\ge\dfrac{k+2}{2},\;\;k\ge
0.
\end{equation}
Starting from an arbitrary initial point $x^0$, the fast gradient
algorithm generates a
sequence $\{y^k\}$ of approximate solutions with $y^0=x^0$, and a
sequence $\{x^k\}$ of search points according to the following rule \footnote{The gradient method in
  \cite{Nes83} computes a steplength $\alpha_k$ which can be taken as
  the inverse of the Lipschitz constant $L$ for the gradient of the
  objective when it is available. Since $F_\mu$ is $\mu$-smooth, the Lipschitz constant is 
$\mu$. We avoid also the evaluations of $F_\mu$-values  
which would be necessary if a steplength has to be computed.} 
(when applied to $F_\mu$),
\begin{equation}\label{nagd}
y^{k+1}=x^k-\dfrac{1}{\mu}\nabla F_{\mu}(x^k)=p_\mu(x^k),\; 
x^{k+1}=y^{k+1}+\alpha_k(y^{k+1}-y^k), 
\;\alpha_k=\lambda_{k+1}^{-1}(\lambda_k-1).
\end{equation}
The above scheme is in fact a special result of Nesterov's key idea of forming
estimating sequences to devise optimal first-order methods for smooth
optimization \cite{Nes04}. The improvement over the gradient descent
relies on the introduction of the \emph{momentum term} $y^{k+1}-y^k$
as well as the particular coefficient from the sequence
$\{\lambda_k\}$. By considering \emph{one} smooth optimization problem 
\eqref{MYpb} on which
Nesterov's scheme is applied, our goal is not to have to 
 tune the proximity parameter $\mu$,  
looking for acceleration
through the momentum term only.     
In terms of the minimization problem \eqref{pb}, it certainly 
makes sense to consider a varying parameter $\mu$, say
$y^{k+1}=p_{\mu_k}(x^k)$. 
By doing so, 
the resulting scheme writes
\begin{equation}\label{nagd-ipa}
  y^{k+1}=\arg\min\limits_{x\in\rit^n}\left\{f(x)+\dfrac{\mu_k}{2}\|x-x^k\|^2\right\},\; 
  x^{k+1}=y^{k+1}+\alpha_k(y^{k+1}-y^k), 
\end{equation}
and can be cast into 
the recent class of so-called \emph{inertial proximal methods}, the origins of which go back to
\cite{AlAt01}. There is a rich literature devoted to this class of methods, see for
instance the recent papers \cite{ACR18,ACR19}. They stem from the use of an implicit
discretization of a differential system of second-order in time and
give interesting insight into Nesterov's scheme \cite{SBC14}. Applied
to the original problem \eqref{pb}, an
iteration of the inertial proximal algorithm with parameters
$\alpha_k\ge 0$ and  $\tau_k=\mu_k^{-1}>0$ is given (with our notations) by
\eqref{nagd-ipa} but with a \emph{more general 
nonnegative sequence} $\{\alpha_k\}$ of extrapolation
coefficients (including $\{\lambda_{k+1}^{-1}(\lambda_k-1)\}$) that capture the inertial effect of the differential
system. The sequence $\{\tau_k\}$ is interpreted as a sequence of proximal parameters
taking into account the temporal scale effects of the system. 
In fact, the main proximal algorithm proposed by G\"uler for \eqref{pb} can also be written as an
inertial proximal algorithm with some appropriate parameters, 
see \cite{ACR18,ACR19}. 
Although taking insipration from Nesterov's method in
\cite{Nes83}, the second algorithm proposed by G\"uler for \eqref{pb} in Section~6 of \cite{Gul92}
  is different from the scheme \eqref{nagd}, in that 
the update of $x^{k+1}$ involves $x^k$ as follows,
\begin{equation}\label{guler-update}
x^{k+1}=y^{k+1}+\dfrac{\lambda_k-1}{\lambda_{k+1}}(y^{k+1}-y^k)+\dfrac{\lambda_k}{\lambda_{k+1}}(y^{k+1}-x^k).
\end{equation}
This is the same rule for the update of the sequence $\{x^k\}$ proposed in the recent work by Kim and Fessler for
their proposed optimized gradient method OGM1, see [\cite{KiFe16}, page 99]. The
authors seem not to know the work of G\" uler \cite{Gul92}, it is
not referenced in their paper. 
G\" uler proposed this rule intuitively with no
explanation, while in \cite{KiFe16}, it is shown that it corresponds to an
optimal choice of parameters obtained through a relaxed \emph{performance estimation
  problem} introduced by Drori and Teboulle in \cite{DrTe14} to
optimize first-order algorithms. 
We will consider the update \eqref{guler-update} as well for a
second algorithm 
 through the following general rule,
\begin{equation}\label{one-rule}
x^{k+1}=y^{k+1}+\alpha_k(y^{k+1}-y^k)+\beta_k(y^{k+1}-x^k),
\end{equation}
where $\alpha_k=\lambda_{k+1}^{-1}(\lambda_k-1)$ is Nesterov's extrapolation coefficient given in
\eqref{nagd} and $\{\beta_k\}_{k\ge 0}$ is one of the two sequences~:  
\begin{itemize}
\item[\personalitem] $\beta_k=0,\;k\ge 0$ (then \eqref{one-rule}
  reduces to the update of $x^{k+1}$ in \eqref{nagd}) or, 
\item[\personalitem] $\beta_k=\lambda_k\lambda_{k+1}^{-1},\;k\ge 0$
  (to get \eqref{guler-update}).  
\end{itemize}
We are now ready to propose a conceptual fast algorithm for 
the minimization of the smooth function $F_\mu$ and consequently for solving \eqref{pb}. 
\begin{algorithm}\label{fast-ppa}
\begin{enumerate}
\item[] 
\hspace*{-0.5cm}{\tt Fast Proximal Point Algorithm (FPPA)}
\item[0.]\label{fast-ppa-1a} Choose $x^0=y^0\in\rit^n$ and the sequence
  $\{\beta_k\}_{k\ge 0}$. Set $k=0$.  
\item\label{fast-ppa-1b}Set $y^{k+1}=p_\mu(x^k)$.
\item\label{fast-ppa-1c} If $y^{k+1}=x^k$ stop: $x^k$ solves \eqref{pb}.
\item\label{fast-ppa-1d} 
Set $x^{k+1} = y^{k+1}+\alpha_k(y^{k+1}-y^k) +\beta_k(y^{k+1}-x^k)$.  
\item Increase $k$ by 1 and loop to Step \ref{fast-ppa-1b}.
\end{enumerate}
\end{algorithm}
We refer to the two algorithms depending on the choice of $\beta_k$
respectively by FPPA1 and FPPA2. 
The main difference of FPPA with PPA is that the proximal point is not
computed at the previous iterate but rather at a specific linear
combination  of the two previous proximal points for FPPA1, and a
second momentum term $y^{k+1}-x^k$ with the coefficient $\lambda_k\lambda_{k+1}^{-1}$ for
FPPA2. 
\subsection{An inexact first-order oracle for $F_\mu$}
Now, FPPA is not implementable as is since obtaining the exact proximal point $p_\mu(x)$ for any given
$x\in\rit^n$ is as difficult as solving the original problem
\eqref{pb}.  Hopefuly, bundle methods offer a practical mean to
compute a proximal point approximately as follows, see [\cite{HUL93},
Section XV.4.3] and the recent survey \cite{Fra18} on bundle methods. 
Assume a \emph{first-order (exact) oracle} for $f$ is available, that given, $z\in\rit^n$ computes
$f(z)$ and a subgradient $g(z)\in\partial f(z)$. At a given step $j$, after a given number
of calls to the oracle at different points $z^i,\;i=1,\ldots,$
with $g^i\in\partial f(z^i)$, we can form the so-called bundle
$\bundleset_j=\left\{(z^i,f(z^i),g^i)\right\}$ and built the
following approximation function of $f$ defined by,
$$
\fmodel_{\bundleset_j}(x)=\max\left\{f(z^i)+\langle g^i,\; x-z^i\rangle:\;(z^i,f(z^i),g^i)\in\bundleset_j\right\}.
$$
This lower 
approximation ($\fmodel_{\bundleset_j}\le f$) replaces $f$ in \eqref{MY} to yield the following
quadratic problem 
\begin{equation}\label{qp-subpb}
F_{\mu,\bundleset_j}(x)=\min\limits_{z\in\rit^n}\left\{\fmodel_{\bundleset_j}(z)+\dfrac{\mu}{2}\|z-x\|^2\right\}.
\end{equation}
whose solution $z^{j}$ tends to $p_{\mu}(x)$ as the bundle grows, see \cite{Fuk84}. In practice, $z^j$ is considered as an approximation
of $p_\mu(x)$  when the
following condition is met [\cite{HUL93}, Chapter XV],
\begin{equation}\label{stop-bundle}
f(z^{j})-\fmodel_{\bundleset_j}(z^{j})\le\varepsilon,\;\;\varepsilon >0.
\end{equation}
In proximal bundle methods, $\varepsilon$ in \eqref{stop-bundle} is 
usually taken as $\varepsilon=(1-\sigma)\lbrack f(x)-\fmodel_{\bundleset_j}(z^{j})\rbrack$ for some $0<\sigma<1$, in
which case the condition writes
$f(z^{j})\le f(x)-\sigma\lbrack f(x)-\fmodel_{\bundleset_j}(z^{j})\rbrack$,  
resulting in a decrease of the objective function $f$ from
$x$ to $z^j$.  
We note in passing that the criterion \eqref{stop-bundle} used here to
identify an approximate proximal point is clearly much easier to
check than \eqref{rock-cond}. As pointed out in \cite{Fuk84}, it does not imply those of
\cite{Roc76}, in particular \eqref{rock-cond} used in \cite{Gul92}.

An algorithm to compute an approximation $\widecheck 
p_{\mu}(x)$ of the proximal point of a given $x\in\rit^n$ with a  
tolerance $\varepsilon$ is as
follows. We consider it as \emph{the (inexact) 
first-order oracle} for $F_\mu$. 
\begin{algorithm}\label{ppo}
\begin{enumerate}
\item[] 
\hspace*{-0.5cm}{\tt Approximate Proximal Point Oracle (APPO)  at $x\in\rit^n$}
\item[0.]\label{ppo0} Initialize the bundle $\bundleset_j,\; j=1$.  
\item\label{ppo1} Compute the solution 
$z^{j}$ of \eqref{qp-subpb}
\item\label{ppo2} If \eqref{stop-bundle} holds, set 
$\widecheck p_{\mu}(x)=z^{j}$ and exit.
\item\label{ppo3} Compute $f(z^j),\; g^j\in\partial f(z^j)$ and
  incorporate $(z^j,f(z^j),g^j)$ to the bundle Increase $j$ by 1 and loop to Step \ref{ppo1}.
\end{enumerate}
\end{algorithm}
Efficient algorithms have been proposed by  
Frangioni~\cite{Fra96} and Kiwiel~\cite{Kiw94} for solving the special quadratic problem
\eqref{qp-subpb}. We review some basic results of the sequence
generated by APPO useful for our subsequent analysis. Let us introduce the
functions
$$
F_x(z)=f(z)+\dfrac{\mu}{2}\|z-x\|^2\;\;\text{and}\;\; \widecheck{F}_{x,j}(z)=\fmodel_{\bundleset_j}(z)+\dfrac{\mu}{2}\|z-x\|^2.
$$
By definition,
$F_\mu(x)=\min\limits_{z\in\rit^n}F_x(z)=F_x(p_\mu(x))$.
The properties of 
the sequence $\{z^j\}$ generated by the iterative procedure APPO can be found in  
[\cite{Fuk84}, Proposition~3], namely the following
\begin{equation}\label{fuk3}
\widecheck{F}_{x,j}(z^j)\le \widecheck{F}_{x,j+1}(z^{j+1})\le F_\mu(x).
\end{equation}
As the bundle $\bundleset_j$ grows, $f(z^j)$ and $\fmodel_{\bundleset_j}(z^j)$ get closer to each
other i.e. $\lim\limits_{j\to\infty}\lbrack
f(z^j)-\fmodel_{\bundleset_j}(z^j)\rbrack\to 0$. 
The condition \eqref{stop-bundle} is  satisfied for large $j$ when $z^j$ becomes close to $p_\mu(x)$, justifying the
fact that we consider $z^j$ as an approximate proximal point of $x$
when \eqref{stop-bundle} occurs. 
APPO then provides an approximate gradient as $\mu
(x-\widecheck p_{\mu}(x))$ and an approximate function value as 
$F_x(\widecheck p_{\mu}(x))$ since at stop we get 
\begin{equation}\label{Fbound}
F_\mu (x)\le  F_x(\widecheck p_{\mu}(x)) 
\le F_\mu (x)+\varepsilon.
\end{equation}
Indeed, clearly $F_\mu (x)\le  F_x(\widecheck p_{\mu}(x))$. Next,
$$
F_x(\widecheck p_{\mu}(x))  = f(\widecheck
p_{\mu}(x))+\dfrac{\mu}{2}\|\widecheck p_{\mu}(x)-x\|^2
\stackrel{\eqref{stop-bundle}}{\le} \fmodel_{\bundleset_j}(\widecheck p_{\mu}(x))+\dfrac{\mu}{2}\|\widecheck
 p_{\mu}(x)-x\|^2 +\varepsilon =  \widecheck{F}_{x,j}(\widecheck p_{\mu}(x))+\varepsilon.
 $$
We then get \eqref{Fbound}  from the fact that $\widecheck{F}_{x,j}(\widecheck p_{\mu}(x))\le F_\mu(x)$,  see \eqref{fuk3}.
The necessary and sufficient optimality condition for the quadratic
problem~\eqref{qp-subpb} at the stop of APPO (with the bundle set
$\bundleset_j$) writes
$0\in\partial\fmodel_{\bundleset_j}(\widecheck p_{\mu}(x))-\mu(x-\widecheck p_{\mu}(x))$. 
Hence, for any $z\in\rit^n$, we have
$$
f(z)\ge \fmodel_{\bundleset_j}(z) 
\ge 
f(\widecheck p_{\mu}(x))+\langle \mu(x-\widecheck p_{\mu}(x)), z-\widecheck
      p_{\mu}(x)\rangle-\lbrack f(\widecheck p_{\mu}(x))-\fmodel_{\bundleset_j}(\widecheck
      p_{\mu}(x))\rbrack.
$$
and from \eqref{stop-bundle},
$f(z)\ge f(\widecheck p_{\mu}(x))+\langle \mu(x-\widecheck p_{\mu}(x)), z-\widecheck p_{\mu}(x)\rangle-\varepsilon$. 
In other words, 
\begin{equation}\label{subdif}
\mu(x-\widecheck p_{\mu}(x))\in\partial_{\varepsilon} f (\widecheck p_{\mu}(x)).
\end{equation}
It is worth mentioning that APPO is not an inexact first-order oracle in the sense
of \cite{DGN14}. It is also different from the procedure given in Section~3.3 for computing approximate
solutions for the Moreau-Yosida
regularization.  
The inexact oracle for $F_\mu$ proposed in \cite{DGN14} computes a
pair $(F_{\mu,\delta}(x),g_\delta(x))$ which satisfies the following two
inequalities within a tolerance $\delta\ge 0$~: 
$$
0\le F_\mu(z)-\left(F_{\mu,\delta}(x)+\langle
  g_\delta(x),\;z-x\rangle\right)\le\frac{\mu}{2}\|z-x\|^2+\delta,\;\;x,z\in\rit^n, 
$$
which are relaxations of the inequalities
$$
0\le F_\mu(z)-\left(F_\mu(x)+\langle \nabla F_\mu(x),\;z-x\rangle\right)\le\frac{\mu}{2}\|z-x\|^2,\;\;x,z\in\rit^n,
$$
which result from the fact that $F_\mu$ has
Lipschitz continuous gradient.
The provided pair has the following properties. $F_{\mu,\delta}(x)$
is a lower approximation of $F_\mu(x)$ in the following sense
$F_{\mu,\delta}(x)\le F_\mu(x)\le F_{\mu,\delta}(x)+\delta$, 
while $g_\delta(x)$ is a $\delta$-subgradient
of $F_\mu$ at $x$~i.e. $F_\mu(z)\ge F_\mu(x)+\langle g_\delta(x),\;z-x\rangle-\delta,\;z\in\rit^n$. 
Even setting $\varepsilon=\delta$, these features are different from what we have with \eqref{Fbound} and
\eqref{subdif} which are satisfied by the output $(F_x(\widecheck
p_{\mu}(x)),\mu(x-\widecheck p_{\mu}(x)))$ from APPO. Several papers e.g. \cite{SRB11,VSBV13} have been devoted to the study of
errors (in different ways as in the present work) in accelerated proximal gradient methods proposed for the case $f$ is of
the form $f=g+h$ where $g$ and $h$ are convex but $h$ is
differentiable, taking advantage of this structure.  

\subsection{Fast proximal bundle algorithms}
An implementable version of Algorithm \ref{fast-ppa} is obtained by
using APPO for the approximate computation of $p_\mu(x^k)$ in  its step~\ref{fast-ppa-1b}. It
is described as follows. 
\begin{algorithm}\label{fast-pba}
\begin{enumerate}
\item[] 
\hspace*{-0.5cm}{\tt Fast Proximal Bundle Algorithm (FPBA)}
\item[0.]\label{fast-pba-1a} Choose $x^0=y^0\in\rit^n$ and the
  sequence $\{\beta_k\}_{k\ge 0}$. Define the sequence
  $\{\varepsilon_k\}_{k\ge 0}$. 
Set $k=0$.  
\item\label{fast-pba-1b}Call  APPO at $x=x^k$ with $\varepsilon
  =\varepsilon_k$  and set $y^{k+1}=\widecheck  p_\mu(x^k)$.
\item\label{fast-pba-1c} 
Set $x^{k+1} = y^{k+1}+\alpha_k(y^{k+1}-y^k) +\beta_k(y^{k+1}-x^k)$.
\item Increase $k$ by 1 and loop to Step \ref{fast-pba-1b}.
\end{enumerate}
\end{algorithm}
As for FPPA, we refer the two versions of FPBA according to the choice
of the sequence $\{\beta_k\}$ to FPBA1 and FPBA2 respectively. The
latter can be viewed as an implementable version of G\"uler second
algorithm if a fixed parameter $\mu$ is considered (in \cite{Gul92}, it is allowed
to depend on $k$). 
The work performed at a previous call to APPO can be
exploited in the initialization of the bundle at a next call. 
The algorithm FPBA is presented below in the usual description of proximal bundle algorithms. It
involves inner iterations (corresponding to the so-called \emph{null steps}) implementing APPO and outer
iterations (\emph{descent} or \emph{serious steps}) for the generation of the sequences
$\{y^k\}$ and $\{x^k\}$. In this form, the main difference with the
standard proximal bundle algorithm lies in the \emph{stability
center} $x^k$ which is usually taken from (in our notations) the sequence $\{z^j\}_{j\le
  k}$ even if this is not necessary to get convergence, see
\cite{AFFG13,Fra18}. Here, $x^k$ is obtained from a fast gradient iteration.  Also, to the contrary of a classical proximal bundle algorithm, 
at each serious step there is no guarantee of decrease in the objective function value between two successive approximate solutions $y^k$ and $y^{k+1}$
(in \cite{AFFG13} a serious step
does not correspond to a decrease in the objective value as well). 
\begin{algorithm}\label{newprox}
\begin{enumerate}
\item[] 
\item[0.]\label{fpcpa-stepi}  Choose an initial point $x^0\in \rit^n$ and the sequence $\{\beta_k\}_{k\ge 0}$. Define the sequence
  $\{\varepsilon_k\}_{k\ge 0}$. 
Set $y^0=x^0$, 
$k=0$ and $\lambda_0=1$.
\item\label{fpcpa-step0} Set $z^0=x^k$ and set $j=0$. Compute $f(z^j),\;
  g^j\in\partial f(z^j)$ and initialize $\mathcal B_j$.
\item\label{fpcpa-step1}  If $g^{j}=0$, terminate: $z^j$ solves \eqref{pb}.
\item\label{fpcpa-step2} Get the solution $z^{j+1}$ of the quadratic problem
$$
\min\limits_{z\in\rit^n}\left\{\fmodel_{\bundleset_j}(z)+\dfrac{\mu}{2}\|z-x^k\|^2\right\},
$$
\item\label{fpcpa-step6} Compute $f(z^{j+1})$ and $g^{j+1}\in\partial
  f(z^{j+1})$.
  
If $f(z^{j+1})-\fmodel_{\bundleset_j}(z^{j+1})\le \varepsilon_k$ then ($\widecheck p_\mu(x^k)$ is computed)
\begin{itemize}
\item[\personalitem] Set 
$$
\lambda_{k+1}=\dfrac{1+\sqrt{1+4\lambda_{k}^2}}{2}
$$
\item[\personalitem] Set $y^{k+1}=z^{j+1},\;\;x^{k+1}
  =y^{k+1}+\alpha_k(y^{k+1}-y^k)+\beta_k(y^{k+1}-x^k)$.
\item[\personalitem] Set $\bundleset_k=\bundleset_j,\;k=k+1$ and go to Step \ref{fpcpa-step0}.
\end{itemize}
Otherwise, set $\bundleset_{j+1}=\bundleset_j\cup\{(z^{j+1},f(z^{j+1}),
  g^{j+1})\}$, increase $j$ by 1 and loop to Step~\ref{fpcpa-step1}.
\end{enumerate}
\end{algorithm}
The \emph{subgradient selection} or
\emph{subgradient aggregation} techniques may be used to maintain the size of the bundle
reasonable, see for instance \cite{Kiw90}. 
Using the fact that $y^{k+1}$ solves the quadratic problem in Step~\ref{fpcpa-step2}, we have
\begin{equation}\label{positive-diff}
  \fmodel_{\bundleset_k}(y^{k+1}) \le
  \fmodel_{\bundleset_k}(y^{k+1})+\dfrac{\mu}{2}\|y^{k+1}-x^k\|^2 \stackrel{\eqref{fuk3}}{\le} F_\mu(x^k) \le   f(x^k).
\end{equation}
If it happens that
\begin{equation}\label{x-optimal}
  f(x^k)-\fmodel_{\bundleset_k}(y^{k+1}) \le\eta,
\end{equation}
for some $\eta\ge 0$, then
\begin{equation}\label{x-optimality}
F_\mu(x^k)\le f(x^k) \stackrel{\eqref{x-optimal}}{\le} \fmodel_{\bundleset_{k}}(y^{k+1})+\eta \stackrel{\eqref{positive-diff} }{\le} F_\mu(x^k)+\eta.
\end{equation}
So, when \eqref{x-optimal} holds and $\eta$ is sufficiently small, we may conclude that $f(x^k)\approx
F_\mu(x^k)$ and then $x^k$ solves approximately \eqref{pb} according to 
Proposition~\ref{MYprop}. Another consequence of  \eqref{x-optimal} is the following relation,
$$
f(x^k)\le f(z)+\sqrt{2\eta\mu}\|z-x^k\|+\eta\;\;\text{for
  all}\;\;z\in\rit^n,
$$
which is used sometimes to show that \eqref{x-optimal} implies the
(approximate) optimality of $x^k$ 
if $\eta$ is small, see for instance
\cite{Fuk84,ChFu99}. However, 
for large $\mu$, $\sqrt{2\eta\mu}$ may not be negligeable even
if $\eta$ is very small. The relation
\eqref{x-optimal} is enough on its own as shown by
\eqref{x-optimality}. Note also that \eqref{x-optimal} implies
$$
\fmodel_{\bundleset_k}(y^{k+1}) \le
f(x^k)\le f(y^{k+1})+\eta\le \fmodel_{\bundleset_k}(y^{k+1})+\varepsilon_k+\eta,
$$
where the last inequality comes from the definition of $y^{k+1}$. 
Hence, if $\varepsilon_k$ is small as well, $y^{k+1}$ could also be
considered as an approximate solution. 

We finally observe that, if we discard the momentum (i.e. $\alpha_k=\beta_k=0$), Algorithm~\ref{newprox} becomes a proximal
bundle algorithm with a fixed penalty parameter. The present approach can be extended to 
convex optimization methods related to proximal point algorithms such as those
proposed in \cite{EcBe92,MOP95,Spin85}.

\section{Convergence analysis}\label{sec-conv}
In this section, we consider the global convergence and rate of convergence
of FPBA in its two variants. Let $x^*\in\X^*$ and denote $R=\|x^0-x^*\|(=\|y^0-x^*\|)$. The application of 
[\cite{Nes83}, Theorem 1] in a straightforward manner to FPPA1 (which
aims at solving \eqref{MYpb}) gives the
following convergence estimate
$$
F_\mu(y^{k})-f^*\le \dfrac{4\mu R^2}{(k+2)^2},
$$
where $x^*$ is any optimal solution of \eqref{MYpb}
(and so that of \eqref{pb}) and $f^*=f(x^*)$.
This estimate has been improved by Beck and Teboulle 
in \cite{BeTe09} to
\begin{equation}\label{BT-compl}
F_\mu(y^{k})-f^*\le \dfrac{2\mu R^2}{(k+1)^2}.
\end{equation}
The second algorithm proposed by G\"uler [\cite{Gul92}, Section 6]
applies to \eqref{pb} and requires  exact proximal points. It uses the
relation \eqref{guler-update} with proximity
parameters depending on $k$ and decreasing i.e. 
$\mu_{k+1}\le\mu_k,\;k\ge 0$ starting from some $\mu_0>0$. The
complexity estimate of this algorithm involves the original function as,
\begin{equation}\label{Gubound} 
f(y^k)-f^*\le\dfrac{\mu_0 R^2}{(k+1)^2}.
\end{equation}
Since $F_\mu\le f$ for any $\mu>0$, the same
bound holds for $F_{\mu_0}(y^k)-f^*$, so 
this result improves over \eqref{BT-compl}. The bound obtained in
\cite{KiFe16} improves slighlty on \eqref{Gubound} since in their complexity estimate, 
 $(k+1)(k+1+\sqrt{2})$ replaces $(k+1)^2$ in the r.h.s. of \eqref{Gubound}.

Those bounds do not  apply to FPBA since, to the contrary of FPPA, it
uses approximate  proximal points. We now give 
the convergence results for the two versions of FPBA. For a given iterate $y^k$, let $\delta_k=f(y^k)-f^*$ 
The complexity estimate of the algorithms hinges on a 
lower bound on $\lambda_{k-1}^2\delta_k-\lambda_k^2\delta_{k+1}$. 
We start by giving a common one  
for the sequences $\{y^k\}$ generated by FPBA1 and FPBA2.
\begin{lemma}\label{lem4th1and2}
Assume that the sequence $\{(x^k,y^k)\}$ is generated by FPBA. Then,
$$
\lambda_{k-1}^2\delta_k-\lambda_k^2\delta_{k+1}\ge \mu\langle u^k,v^k\rangle-\lambda_k^2\varepsilon_k, 
$$
where $u^k=\lambda_k(y^{k+1}-x^k)$ and $v^k=\lambda_k(y^{k+1}-y^k)+y^k-x^*$.
\end{lemma}
\proof Using \eqref{subdif} with $x=x^k$ and $\varepsilon=\varepsilon_k$, we have for any $x\in\rit^n$ 
$$
f(x) \ge  f(y^{k+1})+\mu\langle
x^k-y^{k+1},\;x-y^{k+1}\rangle-\varepsilon_k.
$$
We use this inequality with $x=y^k$ and $x=x^*\in\X^*$ to get  respectively
\begin{equation}\label{sgrad-ineq}
\left\{
\begin{array}{l}
f(y^k)-f(y^{k+1})\ge\mu\langle
x^k-y^{k+1},\;y^k-y^{k+1}\rangle-\varepsilon_k,\\
\\[-2mm]
f(x^*)-f(y^{k+1})\ge \mu\langle
x^k-y^{k+1},\;x^*-y^{k+1}\rangle-\varepsilon_k.
\end{array}
\right.
\end{equation}
We proceed as in \cite{Bub15},  
multiplying the first inequality of \eqref{sgrad-ineq} by $\lambda_k-1$ and adding the
result to the second inequality to get
$$
(\lambda_k-1)\delta_k-\lambda_k\delta_{k+1}\ge  \mu\langle
x^k-y^{k+1},\lambda_k(y^k-y^{k+1})+x^*-y^k\rangle-\lambda_k\varepsilon_k.
$$ 
Now, multiplying this inequality by $\lambda_k$, using the relation
\eqref{lrel}, we obtain
$$
\lambda_{k-1}^2\delta_k-\lambda_k^2\delta_{k+1}\ge \mu\langle\lambda_k(x^k-y^{k+1}),\lambda_k(y^k-y^{k+1})+x^*-y^k\rangle-\lambda_k^2\varepsilon_k
 $$
 \endprf
 \noindent We have
 $$
 u^k+v^k=\lambda_k(2y^{k+1}-y^k-x^k)+y^k-x^*,\;\;u^k-v^k=\lambda_k(y^k-x^k)-y^k+x^*.
 $$
 The polarization identity writes
\begin{equation}\label{polar}
\langle u,v\rangle=\frac{1}{4}(\|u+v\|^2-\|u-v\|^2),\;\;u,v\in\rit^n.
\end{equation}
By developping the first term, we have
$$
\langle u,v\rangle  = \frac{1}{4}(\|u+v\|^2-\|u-v\|^2) = \frac{1}{4}(\|u\|^2+2\langle u,v\rangle+\|v\|^2-\|u-v\|^2),
$$
and recover the parallelogram law,
\begin{equation}\label{p-law}
\langle u,v\rangle  = \frac{1}{2}(\|u\|^2+\|v\|^2-\|u-v\|^2). 
\end{equation}
Now, let $w^k=v^k-u^k=\lambda_k(x^k-y^k)+y^k-x^*,\;k\ge 0$. Then,
$$
 w^{k+1}=\lambda_{k+1}(x^{k+1}-y^{k+1})+y^{k+1}-x^*,\;k\ge 0.
$$
According to the updating rules of the proximal point, we have
$$
\lambda_{k+1}(x^{k+1}-y^{k+1})=
\left\{
\begin{array}{l}
  \lambda_k(y^{k+1}-y^k)+y^k-y^{k+1}\;\;\text{ for FPBA1}\\
  \\[-2mm]
\lambda_k(2y^{k+1}-y^k-x^k)+y^k-y^{k+1}\;\;\text{ for FPBA2}.
\end{array}
\right.
$$
Hence, depending on the rule used, $w^{k+1}$ takes another form,
\begin{equation}\label{w-update}
w^{k+1}=
\left\{
\begin{array}{l}
  \lambda_k(y^{k+1}-y^k)+y^k-x^*=v^k\;\;\text{ for FPBA1}\\
  \\[-2mm]
\lambda_k(2y^{k+1}-y^k-x^k)+y^k-x^*=u^k+v^k\;\;\text{ for FPBA2}.
\end{array}
\right.
\end{equation}
Based on this and the common lower bound given in Lemma~\ref{lem4th1and2},
we derive other lower bounds for FPBA1 and FPBA2 involving only the sequence $\{w^k\}_{k\ge 0}$.
\begin{lemma}\label{th1-lem}
Assume that the sequence $\{(x^k,y^k)\}$ is generated by FPBA1. Then
\begin{equation}\label{zeta-eq1}
\lambda_{k-1}^2\delta_k-\lambda_k^2\delta_{k+1}\ge \dfrac{\mu}{2}\|w^{k+1}\|^2-\dfrac{\mu}{2}\|w^k\|^2
-\lambda_k^2\varepsilon_k.
\end{equation}
\end{lemma}
\proof Using \eqref{p-law}, we get
$$
\langle u^k,v^k\rangle=\frac{1}{2}(\|u^k\|^2+\|v^k\|^2-\|u^k-v^k\|^2)\ge \frac{1}{2}(\|v^k\|^2-\|u^k-v^k\|^2).
$$
Therefore,
$$
 \lambda_{k-1}^2\delta_k-\lambda_k^2\delta_{k+1} \ge  \mu\langle u^k,v^k\rangle-\lambda_k^2\varepsilon_k\ge 
\dfrac{\mu}{2}\left[\|v^k\|^2-\|u^k-v^k\|^2\right]-\lambda_k^2\varepsilon_k.
$$
But $u^k-v^k=-w^k$ and for FPBA1, we have $w^{k+1}=v^k$ (see \eqref{w-update}). This gives the desired result.

\endprf
\noindent An analogue result for FPBA2 is as follows. 
\begin{lemma}\label{th2-lem}
 Assume that the sequence $\{(x^k,y^k)\}$ is generated by FPBA2. Then
\begin{equation}\label{zeta-eq2}
\lambda_{k-1}^2\delta_k-\lambda_k^2\delta_{k+1}\ge 
  \dfrac{\mu}{4}\|w^{k+1}\|^2-\dfrac{\mu}{4}\|w^{k}\|^2
-\lambda_k^2\varepsilon_k.
\end{equation}
\end{lemma}
\proof Noting that for FPBA2, $w^{k+1}=u^k+v^k$, we get from \eqref{polar},   
$$
\langle u^k,v^k\rangle=\frac{1}{4}\left(\|w^{k+1}\|^2-\|w^{k}\|^2\right).
$$
Apply then Lemma~\ref{lem4th1and2}.

\endprf
\noindent We are now ready to give the complexity estimate for FPBA1, using Lemma~\ref{th2-lem}. 
\begin{theorem}\label{fpba1-conv}
The sequence $\{(x^k,y^k)\}$ generated by FPBA1 satisfies the following
bound 
$$
f(y^k)-f^*\le\frac{2\mu R^2}{(k+1)^2}+\frac{1}{\lambda_{k-1}^2}\sum\limits_{i=0}^{k-1}\lambda_i^2\varepsilon_i,\;\;k\ge
1.
$$ 
\end{theorem}
\prf 
Summing the inequalities \eqref{zeta-eq1} for $i=1,\ldots,k-1$, gives (recall
that $\lambda_0=1$)
$$
\lambda_{k-1}^2\delta_k\le
\delta_1+\frac{\mu}{2}\|w^1\|^2+\sum\limits_{i=1}^{k-1}\lambda_i^2\varepsilon_i-\frac{\mu}{2}\|w^k\|^2\le
\delta_1+\frac{\mu}{2}\|w^1\|^2+\sum\limits_{i=1}^{k-1}\lambda_i^2\varepsilon_i.
$$
From the second inequality of \eqref{sgrad-ineq} with $k=0$, we get
$$
\begin{array}{lcl}
  \delta_1 &\le & \mu\langle x^0-y^{1},\;y^1-x^*\rangle+\varepsilon_0\\
 \\[-2mm]
 &= & -\mu\langle x^0-y^{1},\;x^*-y^1\rangle+\varepsilon_0\\
\\[-2mm]
& \stackrel{\eqref{p-law}}{=} & -\dfrac{\mu}{2}\left[\|x^0-y^{1}\|^2+\|y^1-x^*\|-\|x^0-x^*\|^2\right]+\varepsilon_0\\
\\[-2mm]
& \le & -\dfrac{\mu}{2}\|y^1-x^*\|^2+\dfrac{\mu}{2}\|x^0-x^*\|^2
      +\varepsilon_0
\end{array}
$$
Note that $w^1=\lambda_0(y^1-y^0)+y^0-x^*=y^1-x^*$ since $\lambda_0=1$. Therefore
$$
\delta_1\le -\dfrac{\mu}{2}\|w^1\|^2+\dfrac{\mu}{2}\|x^0-x^*\|^2+\varepsilon_0,
$$
and
\begin{equation}\label{th1-bound}
\lambda_{k-1}^2\delta_k\le
\dfrac{\mu}{2}\|x^0-x^*\|^2+\sum\limits_{i=0}^{k-1}\lambda_i^2\varepsilon_i,
\end{equation}
which combined with the fact that $\lambda_{k-1}\ge (k+1)/2$ gives the
desired result.

\endprf

Thanks to a better lower bound obtained in Lemma~\ref{th2-lem} for the sequence generated by FPBA2, its complexity estimate appears better. 
\begin{theorem}\label{fpba2-conv}
The sequence $\{y^k\}$ generated by FPBA2 satisfies the following bound
$$
f(y^k)-f^*\le\frac{\mu R^2}{(k+1)^2}+\frac{1}{\lambda_{k-1}^2}\sum\limits_{i=0}^{k-1}\lambda_i^2\varepsilon_i,\;\;k\ge
1.
$$ 
\end{theorem}
\proof 
As in the proof of Therorem~\ref{fpba1-conv}, we sum the inequalities \eqref{zeta-eq2} for $i=1,\ldots,k-1$ and get
$$
\lambda_{k-1}^2\delta_{k}\le \delta_1
+\frac{\mu}{4}\|w^{1}\|^2+\sum\limits_{i=1}^{k-1}\varepsilon_i\lambda_{i}^2-\frac{\mu}{4}\|w^{k}\|^2\le \delta_1+\frac{\mu}{4}\|w^{1}\|^2+\sum\limits_{i=1}^{k-1}\lambda_{i}^2\varepsilon_i.
$$
We use again the second inequality of \eqref{sgrad-ineq} for $k=0$ to obtain
$$
\begin{array}{lcl}
\delta_1&\le & \mu\langle x^0-y^{1},\; y^{1} -x^*\rangle+\varepsilon_0
  \\
\\[-2mm]
& = &
      \dfrac{\mu}{4}\|x^0-x^*\|^2-\dfrac{\mu}{4}\|x^0-2y^{1}+x^*\|^2+\varepsilon_0\\
\\[-2mm]
& = & \dfrac{\mu}{4}\|x^0-x^*\|^2-\dfrac{\mu}{4}\|w^1\|^2+\varepsilon_0,
\end{array}
$$
noting that 
$w^1=\lambda_0(x^0+y^0-2y^{1})+x^*-y^0=x^0-2y^1+x^*$. Putting together
the above bound on $\delta_1$ 
and the previous inequality, one gets
\begin{equation}\label{th2-bound}
\delta_{k}\le \frac{\mu}{4 \lambda_{k-1}^2}\|x^0-x^*\|^2+\frac{1}{\lambda_{k-1}^2}\sum\limits_{i=0}^{k-1}\lambda_{i}^2\varepsilon_i
\end{equation}
It remains to use the fact that $\lambda_{k-1}\ge (k+1)/2$, see \eqref{lrel2}. 

\endprf

We have fixed a parameter $\mu>0$ and get 
one smooth optimization problem \eqref{MYpb} on which the fast
gradient concept has been applied. In this way, the number of calls to the $F_\mu$-oracle
APPO is optimized. Of course, $\mu$ has an impact in the efficiency of solving the
quadratic subproblems \eqref{qp-subpb} as well as the number of
calls to the first-order oracle for $f$, which is better to be
minimized. There comes the need to adapt $\mu$ at 
each step although this breaks the philosophy of our approach. 
Following a different approach,
the algorithms proposed by 
G\"uler in \cite{Gul92} use proximity parameters depending on $k$
satisfying the condition (with our notations)
\begin{equation}\label{mu-cond}
  \mu_0=\mu\;\text{for some given}\;\mu>0\;
  \text{and}\;\mu_k\le\mu_{k-1},\;k\ge 1.
\end{equation}
In the present setting, it is also possible to use
different parameters under the same condition. In this case, Step~\ref{fast-pba-1b} of FPBA is modified as follows. 
\begin{enumerate}
\item[$1^\prime.$]Call  APPO at $x=x^k$ and $\mu=\mu_k$. Set $y^{k+1}=\widecheck  p_\mu(x^k)$.
\end{enumerate}
\begin{proposition}\label{rmk3.3}
The complexity estimates of Theorems \ref{fpba1-conv} and
\ref{fpba2-conv} hold if instead of a fixed proximity parameter, we
consider a sequence of positive numbers $\{\mu_k\}_{k\ge 0}$
satisfying \eqref{mu-cond}.
\end{proposition}
\prf 
We consider only Theorem~\ref{fpba1-conv} and show that it 
remains valid with the above modification (the proof for Theorem~\ref{fpba2-conv} is similar). It easily seen that Lemma~\ref{lem4th1and2} and 
Lemma~\ref{th1-lem} hold with $\mu_k$ in place of $\mu$. Based on the fact that  $\mu_k\le\mu_{k-1}$, inquality \eqref{zeta-eq1}
 yields
$$
\begin{array}{lcl}
\lambda_{k-1}^2\delta_k-\lambda_k^2\delta_{k+1} &\ge &\dfrac{\mu_k}{2}\|w^{k+1}\|^2-\dfrac{\mu_k}{2}\|w^k\|^2-\lambda_k^2\varepsilon_k\\
  \\[-2mm]
 & \ge & \dfrac{\mu_k}{2}\|w^{k+1}\|^2-\dfrac{\mu_{k-1}}{2}\|w^k\|^2-\lambda_k^2\varepsilon_k
\end{array}
$$ 
Summing these inequalities for $i=1,\ldots,k-1$, yields
$$
\begin{array}{lcl}
 \lambda_{k-1}^2\delta_k & \le &\delta_1+\dfrac{\mu_0}{2}\|w^{1}\|^2+\sum\limits_{i=1}^{k-1}\lambda_i^2\varepsilon_i-\dfrac{\mu_{k-1}}{2}\|w^k\|^2\\
 \\[-2mm]
& \le  &\delta_1+\dfrac{\mu}{2}\|w^1\|^2+\sum\limits_{i=1}^{k-1}\lambda_i^2\varepsilon_i\;\;\text{(we use $\mu_0=\mu$)}
\end{array}
$$
In the present context, \eqref{sgrad-ineq} with $k=0$ and
\eqref{p-law} give
$$
\delta_1\le -\dfrac{\mu_0}{2}\|y^1-x^*\|^2+\dfrac{\mu_0}{2}\|x^0-x^*\|^2+\varepsilon_0 =  -\dfrac{\mu}{2}\|w^1\|^2+\dfrac{\mu R^2}{2}+\varepsilon_0,
$$
and then
$$
\lambda_{k-1}^2\delta_k \le \dfrac{\mu R^2}{2}+\sum\limits_{i=0}^{k-1}\lambda_i^2\varepsilon_i.
$$
\endprf

\begin{remark}\label{rmk3.1lambda}
 Observe that \eqref{th1-bound} and \eqref{th2-bound} remain valid if the sequence $\{\lambda_k\}$
 satisfies the relation
 \begin{equation}\label{gen-l-cond}
 \lambda_k^2-\lambda_{k-1}^2\le \lambda_k,\;k\ge 1,
 \end{equation}
 used in \cite{AC18,ACR18,ACR19} to generalize the
 extrapolation coefficients $\alpha_k=\lambda_{k+1}^{-1}(\lambda_k-1)$ for inertial proximal methods.  
 Equality holds in \eqref{gen-l-cond} for Nesterov's
 sequence, cf \eqref{lrel} used in Lemma~\ref{lem4th1and2} (which holds with \eqref{gen-l-cond}). 
\endprf
\end{remark}


\begin{remark}\label{rmk3.2}
 Initially, FPBA intends to solve the minimization problem of
 $F_\mu$. However, the complexity estimates are expressed in terms of $f$-values. If we discard the errors 
in these complexity estimates, we recover the known ones given at the begining of this section for FPPA1 and FPPA2 using exact
proximal points.  
One cannot draw a conclusion of the superiority of a scheme
to the other from the above complexity estimates. These worst-case
convergence bounds are the ones 
we were able to establish. 
We cannot exclude that it is possible to
get tighter bounds. \endprf
\end{remark}

\begin{remark}\label{bundle-complexity}
  Complexity estimates for classical proximal bundle methods have been established 
  requiring $O(\varepsilon^{-2})$ outer iterations and $O(\varepsilon^{-3})$ iterations while taking into account the number of inner iterations. 
 For proximal level bundle methods, the complexity estimate is $O(\varepsilon^{-2})$. See for instance \cite{Fra18,Kiw00}. 
  A subsequent work is needed to include inner iterations in the complexity analysis of FPBA. 
\end{remark}

\section{Error accumulation}\label{sec-err}
It is pointed out in \cite{DGN14}  that fast first-order methods suffer from
accumulation of errors to the contrary of classical gradient
methods (see also the inexact approach in \cite{Gul92}).
The accumulation of errors at step $k$,
$\vartheta_k=\lambda_{k-1}^{-2}\sum\limits_{i=0}^{k-1}\lambda_{i}^2\varepsilon_i$, 
is identical in both schemes FPBA1 and FPBA2
and similar to that of the fast gradient method with the inexact oracle
proposed in \cite{DGN14}. Since $\lambda_{k-1}\ge (k+1)/2$, we have, 
\begin{equation}\label{acc-error}
\vartheta_k\le\frac{4}{(k+1)^2}\sum\limits_{i=0}^{k-1}\lambda_{i}^2\varepsilon_i.
\end{equation}
\subsection{Error weights}
Let $\omega_{i,k} =\lambda_i^2\lambda_{k-1}^{-2}, \;i=0\ldots,k-1,$ be the weight of the error $\varepsilon_i$ in $\vartheta_k$
(note that it depends on $k$). 
Using the first relation \eqref{lrel}, we have for $i=0\ldots,k-2$,
$$
\omega_{i,k}=\frac{\lambda_i^2}{\lambda_{k-1}^2}=\frac{\lambda_{i+1}^2-\lambda_{i+1}}{\lambda_{k-1}^2}=\omega_{i+1,k}-\frac{\lambda_{i+1}}{\lambda_{k-1}^2}\;\;\text{i.e.}\;\;
\omega_{i+1,k}=\omega_{i,k}+\frac{\lambda_{i+1}}{\lambda_{k-1}^2}.
$$
Hence, $\omega_{i,k}$ increases strictly with $i$ but is bounded by $1$,
$$
0<\omega_{0,k}=\frac{1}{\lambda_{k-1}^2}<\omega_{1,k}<\omega_{2,k}<\ldots 
<\omega_{k-1,k}=1. 
$$
However, $\omega_{i,k}$ decreases with $k$ as
$\lambda_{k-1}$ is increasing. But for a given $k$, we have
$\omega_{i,k}<\omega_{i+1,k}$ for $i=0,\ldots,k-2$ i.e 
the weight increases from $\lambda_{k-1}^{-2}$ to the maximum $\omega_{k-1,k}=1$ (with the weights
in the r.h.s of \eqref{acc-error},
the last ones exceed $1$). With this
observation, one can tolerate large errors in early iterations 
but require smaller and smaller errors in the progress of the algorithms. 
\subsection{Special cases}
\subsubsection{Equal errors}Assume that  $\varepsilon_i=\varepsilon$ for all $i\ge
0$. Based on the first relation in \eqref{lrel2}, we have $\vartheta_k=\theta_k\varepsilon$ where
$$
\theta_k\triangleq\sum\limits_{i=0}^{k-1}\omega_{i,k}=\frac{1}{\lambda_{k-1}^2}\sum\limits_{i=0}^{k-1}\sum\limits_{l=0}^{i}\lambda_{l}=
\frac{1}{\lambda_{k-1}^2}\sum\limits_{i=0}^{k-1}(k-i)\lambda_i= 
k-\frac{1}{\lambda_{k-1}^2}\sum\limits_{i=1}^{k-1}i \lambda_i.
$$
Hence, $\theta_k$ is far away from $k$ and so is the accumulated error $\vartheta_k=\theta_k\varepsilon$ from $k\varepsilon$. But it is
asymptotically divergent with the first terms in the
complexity bounds as it is the case for the fast gradient method of
\cite{DGN14}. Indeed, starting from $\theta_1=1$, $\theta_k$ is
increasing with $k$ as it is shown next. We have
for any $k\ge 1$,
$$
\theta_{k+1}=\frac{1}{\lambda_k^2}\sum\limits_{i=0}^k\lambda_i^2=1+\frac{\lambda_{k-1}^2}{\lambda_k^2}\sum\limits_{i=0}^{k-1}\frac{\lambda_i^2}{\lambda_{k-1}^2}= 
1+(1-\frac{1}{\lambda_k})\theta_k.
$$
Hence, $\theta_{k+1}-\theta_k=1-\lambda_k^{-1}\theta_k$. 
We prove by induction that this difference is
positive. It is true for $k=1$ since $1-\lambda_1^{-1}\theta_1=1-\lambda_1^{-1}>0$.
Assume that  it holds for $k$ i.e. $1-\lambda_k^{-1}\theta_k\ge 0$ and
let show that it holds for $k+1$. We have,
$$
\begin{array}{lcl}
1-\dfrac{1}{\lambda_{k+1}}\theta_{k+1} & =&
                                       1-\dfrac{1}{\lambda_{k+1}}\left\lbrack 1+(1-\dfrac{1}{\lambda_k})\theta_k\right\rbrack\\
\\[-2mm]
& = & 1-\dfrac{1}{\lambda_{k+1}}\theta_k-\dfrac{1}{\lambda_{k+1}}(1-\dfrac{1}{\lambda_k}\theta_k)\\
\\[-2mm]
& \ge &
        1-\dfrac{1}{\lambda_{k}}\theta_k-\dfrac{1}{\lambda_{k+1}}(1-\dfrac{1}{\lambda_k}\theta_k)
        \;\;\text{(as $\lambda_k<\lambda_{k+1}$)}\\
\\[-2mm]
& = & (1-\dfrac{1}{\lambda_{k+1}})(1-\dfrac{1}{\lambda_k}\theta_k)\\
\\[-2mm]
& \ge & 0 \;\;\text{(as $1<\lambda_{k+1}$)}.
\end{array}
$$
The divergence between the two terms in the complexity estimates may
be avoided if $\varepsilon_k=O(k^{-\sigma})$ for some parameter
$\sigma>0$, see \cite{Gul92,DGN14}.

\subsubsection{Step dependent errors}For the case where
$\varepsilon_k$ is different for each step $k$, the sequence
$\{\vartheta_k\}_{k\ge 0}$ satisfies the relation
$$
\vartheta_{k+1}=\varepsilon_k+\frac{\lambda_{k-1}^2}{\lambda_k^2}\sum\limits_{i=0}^{k-1}\frac{\lambda_i^2}{\lambda_{k-1}^2}\varepsilon_i 
\stackrel{\eqref{lrel}}{=}\varepsilon_k+\left(1-\lambda_k^{-1}\right)\vartheta_k.
$$
Note that for $k\ge 1$, we have $1-\lambda_k^{-1}>0$. Since
$\vartheta_{k+1}-\vartheta_k=\varepsilon_k-\lambda_k^{-1}\vartheta_k$, 
the sequence $\{\vartheta_k\}_{k\ge 0}$ may be made decreasing by choosing
$\varepsilon_k\le \lambda_k^{-1}\vartheta_k$ for $k\ge 1$. In this 
case, as $\vartheta_1=\varepsilon_0$, we have
$\vartheta_k\le\varepsilon_0,\;k\ge 1$ and the
complexity estimates of FPBA1 and FPBA2 write respectively
$$
f(y^k)-f^*\le
\frac{2\mu R^2}{(k+1)^2}+\varepsilon_0\;\;\text{and}\;\;
f(y^k)-f^*\le \frac{\mu R^2}{(k+1)^2}+\varepsilon_0.
$$
In other words, there is no accumulation error in this case, 
and $f(y^k)-f^*$ tends asymptotically 
to $\varepsilon_0$.  
In particular, if we set
$\varepsilon_k=\lambda_k^{-1}\vartheta_k,\;k\ge 1$, we
have $\vartheta_k=\varepsilon_0$ for any $k\ge 1$ and therefore
\begin{equation}\label{tol-choice}
\varepsilon_k=\frac{1}{\lambda_k}\varepsilon_0\stackrel{\eqref{lrel2}}{\le}
\frac{2}{k+2}\varepsilon_0.
\end{equation}
If we wish the residual $f(y^k)-f^*$ to reach an accuracy
$\varepsilon$ with FPBA2 for instance, then we
set $\varepsilon_0=\frac{\varepsilon}{2}$ and the number $k$ of steps to perform
should satisfy
$$
\frac{\mu R^2}{(k+1)^2}\le\frac{\varepsilon}{2},
$$
which gives $k\ge R\sqrt{\frac{2\mu}{\varepsilon}}-1$. It easy to
check that for FPBA1, the condition on $k$  is  
$k\ge 2R\sqrt{\frac{\mu}{\varepsilon}}-1$.

The choice \eqref{tol-choice} results in a strictly decreasing errors sequence and the approach FPBA is asymptotically an ``almost exact'' fast gradient method.
It is much interesting to exploit the fact that the weights
of former errors are decreasing to zero as the iterations progress, and then choose the errors in order 
to escape from the bundle mechanism as soon as possible as in classical proximal bundle algorithms.
For instance, it is still possible to use the condition of classical proximal bundle algorithms, 
\begin{equation}\label{bundlestop}
f(z^{j+1})\le f(x^k)-\sigma [f(x^k)-\fmodel_{\bundleset_j}(z^{j+1})],\;0<\sigma<1, 
\end{equation}
and set $y^{k+1}=z^{j+1}$ when it is satisfied, implying $f(y^{k+1})\le f(x^k)$. This would mean setting
$$
\varepsilon_k=(1-\sigma) [f(x^k)-\fmodel_{\bundleset_k}(y^{k+1})],\;\;k\ge 0.
$$

\section{Numerical experiments}\label{sec-num}
In this section, we present some numerical results to provide a first idea about the performance of the proposed
algorithms as compared to some previous proximal algorithms. To this aim, we consider fifteen of the academic test problems
already used 
in \cite{Ouo09}.  
For the proximity parameter, we consider in all our runs  the standard
choice $\mu=1$ which usually suits for  well-scaled problems. 
As the optimal values of the test problems are available, we stop the
algorithms when
$$
f^k_{\text{best}}-f^*\le 10^{-6}(1+|f^k_{\text{best}}|),
$$
where $f^k_{\text{best}}$ is the best function value recorded during
the $k$ steps, or when $\|g^j\|\le 10^{-6}$ for some $j$. Clearly, there is a
need for a practical condition identifying $y^k$ as an approximate solution of \eqref{pb} other than fixing a
number of steps to perform as in \cite{DrTe14,KiFe16} or considering
the r.h.s in the
complexity estimates of Theorems~\ref{fpba1-conv}
and~\ref{fpba2-conv}, which correspond to the worst case performance of the
algorithms. One possibility is to use the following upper approximation of the gap function $\delta(x)$,
$[0\le \delta(x)\le ]\;\;\delta_{\bundleset_k} (x)=f(x)-F_{\mu,\bundleset_k}(x)$,
by checking $\delta_{\bundleset_k} (y^{k+1})\le\varepsilon$ 
for a given precision $\varepsilon>0$, but at the cost of computing $F_{\mu,\bundleset_k}(y ^{k+1})$ at each step $k$. 

We have implemented the algorithms using Python 3.5 and 
Cplex~12.7.1  as the solver of the quadratic problem \eqref{qp-subpb}
which has been reformulated as
\begin{equation}\label{qp-subpb-eq}
\min\left\{w+\dfrac{\mu}{2}\|x-x^k\|^2: \; 
  f(z^i)+\langle g^i,\; x-z^i\rangle\le w,\;\;i\in\bundleset_j,\;w\in\rit,\;\;x\in \rit^n
\right\}
\end{equation}
\begin{table}[htbp]
\footnotesize
\caption{\small Test problems}
\begin{center}
\begin{tabular}{|c|c|c|c|}
\hline
 Problem & Name & $n$ & $f^*$  \\
\hline
1 & CB2 & 2 & 1.952224 \\
\hline
2 & CB3 & 2 & 2  \\
\hline
3 & DEM & 2 & -3  \\
\hline
4 & QL & 2 & 7.2 \\
\hline
5 & LQ & 2 & -$\sqrt{2}$ \\
\hline
6 & Mifflin1 & 2 & -1\\
\hline
7 & Mifflin2 & 2 & -1\\
\hline
8 & Rosen-Suzuki & 4 & -44\\
\hline
9 & Shor & 5 & 22.600162\\
\hline
10 & Maxquad & 10 & -0.841408\\
\hline
11 & Maxq & 20 & 0 \\
\hline
12  & Maxl & 20 & 0 \\
\hline
13  & Goffin & 50 & 0 \\
\hline
14  & MxHilb & 50 & 0 \\
\hline
15  & L1Hilb & 50 & 0 \\
\hline
\end{tabular}
\label{tpb}
\end{center}
\end{table}   
We fix the maximum number of $k$-steps to $250$ in all the runs. 
The results obtained by the two versions of
FPBA are collected on
Table~\ref{tolres} with different values for 
$\varepsilon_0$ in \eqref{tol-choice} whose  r.h.s is taken as $\varepsilon_k$. We reported the number of calls ($\#fg$) to the first-order $f$-oracle for
function and subgradient evaluations at trial solutions $z^j$, the number of steps used by
the algorithms to reach the above stopping criterion
($\#k$). Column $f-f^*$ gives respectively the (absolute) 
difference between the best function value found by the
algorithms at termination and the optimal value. 
The numerical experiments tend to confirm our observation at the end
of Remark~\ref{rmk3.2}. 
At a first glance on  the complexity estimates of Theorems~\ref{fpba1-conv} and~\ref{fpba2-conv}, one would 
expect FPBA2 to outperform FPBA1. We can observe that this is not
 the case since there is no clear superior algorithm among the two 
versions of FPBA, in terms of number of calls to APPO as well as the
number of calls to the first-order oracle for
$f$.  The latter seems to increase with $\varepsilon_0$ for most of the test
problems (we didn't include the results obtained with $\varepsilon_0$ for space limitation).
\begin{table}[htbp]
\scriptsize
\begin{center}
\caption{\small Results with different values of $\varepsilon_0$} 
\begin{tabular}{|c||c|c|c||c|c|c||}
\hline
&\multicolumn{6}{|c||}{$\varepsilon_0 =10^{-1}$} \\
\cline{2-7} \raisebox{1.5ex}[0cm][0cm]{Pb}  & \multicolumn{3}{|c||}{FPBA1} & \multicolumn{3}{|c||}{FPBA2} \\  
\cline{2-7}  & $\#k$ &$\#fg$ & $f-f^*$ & $\#k$ &$\# fg$ & $f-f^*$\\
\hline
1 & 11  & 20 & 1.01E-06 & 15 & 24 &  4.74E-08 \\
\hline
2 & 6 & 13 & 5.66E-09 & 6 & 13 & 5.87E-09\\
\hline
3 & 5 & 10 & 4.37E-11 & 4 & 9 & 1.89E-08 \\
\hline
4 & 9 & 20 & 5.47E-06 & 13 & 24 & 2.70E-06\\
\hline
5 & 2 & 6 & 1.29E-07 & 3 & 8 & 4.31E-08\\
\hline
6 & 11 & 26 &  1.27E-06 & 14 & 28 &  8.60E-07\\
\hline
7 & 15 & 27 & 1.05E-06 & 19 & 31 & 5.63E-08 \\
\hline
8 & 22 & 48 & 2.53E-05& 22 & 50 & 2.87E-05 \\
\hline
9 & 22 & 52 &  1.22E-05 & 23 & 59 & 2.19E-05\\
\hline
10 & 106 & 182 & 8.65E-07 & 140 & 254 & 1.64E-06 \\
\hline
11 & 222 & 491 & 9.26E-07 & 217 & 429 & 9.36E-07\\
\hline
12 & 37 & 77 & 1.22E-08 & 65 & 105 & 2.50E-09 \\
\hline
13 & 11 & 62 &  6.26E-09 & 13 & 64 & 1.50E-08\\
\hline
14 & 206 & 212 & 8.87E-07  & 153 & 160 & 9.99E-07\\
\hline
15 & 72 & 85 &  9.70E-07 & 52 & 65 & 8.75E-07\\
\hline\\
\hline
& \multicolumn{6}{|c||}{$\varepsilon_0 =10^{-3}$} \\
\cline{2-7} \raisebox{1.5ex}[0cm][0cm]{Pb}  & \multicolumn{3}{|c||}{FPBA1} & \multicolumn{3}{|c||}{FPBA2} \\  
\cline{2-7}  & $\#k$ &$\#fg$ & $f-f^*$ & $\#k$ &$\# fg$ & $f-f^*$ \\
\hline
  1 & 8  & 25 & 6.33E-07 & 7 & 26 & 1.99E-06\\
\hline
2 & 3 & 13 & 5.69E-09 & 3 & 13 & 5.57E-09\\
\hline
3 & 4 & 10 & 2.91E-09 & 3 & 9 & 3.36E-06\\
\hline
4 & 11 & 37 & 4.52E-06 & 7 & 30 & 2.96E-06\\
\hline
5 & 1 & 7 &  1.28E-07 & 2 & 8 &  3.80E-08\\
\hline
6 & 11 & 37 & 1.22E-06 & 7 & 31 & 1.87E-06 \\
\hline
7 & 9 & 39 & 1.69E-06 & 13 & 39 & 8.94E-07 \\
\hline
8 & 2 & 57 & 3.60E-05& 9 & 77 & 3.35E-05\\
\hline
9 & 8 & 68 & 1.28E-05 & 8 & 96 & 2.32E-05 \\
\hline
10 & 37 & 187 & 1.30E-06 & 53 & 375 & 1.77E-06\\
\hline
11 & 128 & 967 & 9.09E-07 & 125 & 580 & 7.51E-07\\
\hline
12 & 37 & 78 & 1.22E-08 & 44 & 85 & 1.53E-09\\
\hline
13 & 11 & 62  & 6.26E-09 & 13 & 64 & 1.50E-08\\
\hline
14 & 201 & 230 & 8.97E-07& 147 & 175 & 9.74E-07\\
\hline
15 & 51 & 112 & 9.37E-07& 52 & 109 & 9.49E-07\\
\hline
\hline
\end{tabular}
\label{tolres}
\end{center}
\end{table}
We also experiment the condition \eqref{bundlestop} used in classical proximal bundle algorithms 
with $\sigma=0.5$. 
The results are given on Table~\ref{bmeps} and show that this condition may be a good choice 
as well in the present setting, at the cost of an additional partial
call to the $f$-oracle for the computation of $f(x^k)$. Even by including these calls to count the  number of requests to the oracle,  
escaping from the bundling mechanism as soon as possible may be a winning strategy on some test problems. 
\begin{table}[htbp]
\scriptsize
\begin{center}
\caption{\small Results obtained by FPBA algorithms with the rule \eqref{bundlestop}}
\begin{tabular}{|c||c|c|c||c|c|c||}
\hline
 & \multicolumn{3}{|c||}{FPBA1} & \multicolumn{3}{|c||}{FPBA2} \\  
\cline{2-7} \raisebox{1.5ex}[0cm][0cm]{Pb}  & $\#k$ &$\#fg$ & $f-f^*$ & $\#k$ &$\# fg$ & $f-f^*$\\
  \hline
1 & 10  & 15 & 8.82E-08 & 19 & 22 &  2.82E-06 \\
\hline
2 & 8 & 13 & 5.27E-08 & 9 & 13 & 1.65E-06\\
\hline
3 & 5 & 10 & 4.37E-11 & 5 & 9 & 4.06E-08 \\
\hline
4 & 10 & 16 & 6.17E-06 & 23 & 25 & 3.87E-06\\
\hline
5 & 3 & 7 & 1.29E-07 & 4 & 8 & 3.94E-08\\
\hline
6 & 7 & 24 &  1.17E-06 & 14 & 29 &  8.36E-07\\
\hline
7 & 18 & 30 & 5.88E-07 & 18 & 32 & 8.51E-07 \\
\hline
8 & 13 & 37 & 3.71E-05 & 40 & 51 & 4.31E-05 \\
\hline
9 & 15 & 38 &  2.35E-05 & 48 & 53 & 1.06E-05\\
\hline
10 & 68 & 147 &  1.53E-06 & 257 & 265 & 1.72E-06 \\
\hline
11 & 56 & 202 & 9.67E-07 & 190 & 409 & 6.19E-07\\
\hline
12 & 35 & 75 & 1.22E-08 & 50 & 82 & 1.27E-09 \\
\hline
13 & 14 & 53 &  7.11E-09 & 22 & 64 & 1.39E-08\\
\hline
14 & 203 & 209 & 9.09E-07  &  179 & 184 & 8.76E-07\\
\hline
15 & 57 & 105 &  9.68E-07 & 54 & 57 & 9.37E-07\\
\hline
\end{tabular}
\label{bmeps}
\end{center}
\end{table}
There is certainly a room for improving the practical
efficiency of FPBA by devising practical rules for the management of
the parameter $\mu$ in the lines suggested by Proposition~\ref{rmk3.3} and the popular sequences of the literature on
inertial proximal algorithms. 

Disregarding the way FPBA has been derived, 
other variants of FPBA can
be considered as for proximal bundle algorithm, based on alternative
(equivalent) subproblems of \eqref{qp-subpb}. First, from Proposition~2.2.3 in Chapter~XV of \cite{HUL93},
the exists $\kappa(\mu)>0$ such that any solution of the problem
$$
\min\limits_{z\in\rit^n}\left\{\fmodel_{\bundleset_j}(z):\;\|z-x^k\|^2\le \kappa(\mu)^2\right\},
$$
also solves \eqref{qp-subpb}.  Second, by interpreting $w$ in
\eqref{qp-subpb-eq} (the below 
equivalent reformulation of \eqref{qp-subpb}) as the dualization of a constraint
$w\le l(\mu)$, a \emph{level stabilization} variant of FPBA consits in solving 
$$
\min\left\{\|x-x^k\|^2: \;
f(z^i)+\langle g^i,\; x-z^i\rangle\le l(\mu),\;\;i\in\bundleset_j,\;\; x\in \rit^n
\right\}.
$$
With a suitable choice of $l(\mu)$, the solution of this problem is
that of \eqref{qp-subpb}.
These equivalences are only theoretical as pointed out in \cite{Fra18,HUL93}, 
finding $\kappa(\mu)$ or $l(\mu)$ for a given $\mu$ is not trivial.

Finally, it could be interested to analyze if some improvement on inertial proximal
algorithms may be obtained using a second momentum term, yielding a 
generalized  algorithm
$$
 y^{k+1}=\arg\min\limits_{x\in\rit^n}\left\{f(x)+\dfrac{1}{2\tau_k}\|x-x^k\|^2\right\},
 \; x^{k+1}=y^{k+1}+\alpha_k(y^{k+1}-y^k)+\beta_k(y^{k+1}-x^k).
$$
The  sequence $\{\alpha_k\}$ is general
(including Nesterov's extrapolation coefficients) while $\beta_k$ may be the one we use in this
paper i.e.  $\beta_k=\lambda_k\lambda_{k+1}^{-1}$ since  it is shown in \cite{DrTe14} to correspond to some optimal choice
in first-order algorithms, or any other value that ensures convergence
of the scheme.

\section{Conclusion}\label{sec-con}
We proposed new proximal
bundle algorithms for the minimization of nonsmooth convex
functions, by exploiting fast gradient smooth methods on Moreau-Yosida regularization.
The difference with the proximal bundle algorithm is
the generation of an additional sequence $\{x^k\}$ from which a 
sequence $\{y^k\}$ of proximal points is computed. The computation of $x^k$
is trivial, so the main work is almost the same as in the classical
proximal bundle algorithm. We derive complexity estimates of the
proposed implementable algorithms which suffer from an error accumulation
due to the use of approximate proximal points.

\begin{acknowledgement}
  I'm very grateful to Philippe Mahey for his useful comments on a previous version of the paper, 
\end{acknowledgement}

\end{document}